\documentclass[reqno]{amsart}
\usepackage{graphicx}
\usepackage{pictex}
\usepackage{mathrsfs}
\usepackage{hyperref}
\begin{document}

\newtheorem{thm}{Theorem}[section]
\newtheorem{lem}[thm]{Lemma}
\newtheorem{cor}[thm]{Corollary}
\newtheorem{add}[thm]{Addendum}
\newtheorem{prop}[thm]{Proposition}
\newtheorem{conj}[thm]{Conjecture}
\theoremstyle{definition}
\newtheorem{defn}[thm]{Definition}

\theoremstyle{remark}
\newtheorem{rmk}[thm]{Remark}

%%%%%%%%%%%%%%%%%%%%%%%%%%%%%%%%%%%%%%%%%%%%%%%%%%%%%%%%%

%%\copyrightinfo{2005}{American Mathematical Society}
%%\revertcopyright
%
%\newtheorem{theorem}{Theorem}[section]
%\newtheorem{lemma}[theorem]{Lemma}
%
%\theoremstyle{definition}
%\newtheorem{definition}[theorem]{Definition}
%\newtheorem{example}[theorem]{Example}
%\newtheorem{xca}[theorem]{Exercise}
%
%\theoremstyle{remark}
%\newtheorem{remark}[theorem]{Remark}
%
%\numberwithin{equation}{section}

\newcommand{\OmegaH}{\Omega/\langle H \rangle}
\newcommand{\hatOmegaHstar}{\hat \Omega/\langle H_{\ast}\rangle}
\newcommand{\SurfG}{\Sigma_g}
\newcommand{\TriangG}{T_g}
\newcommand{\TriangGOne}{T_{g,1}}
\newcommand{\ProjG}{\mathcal{P}_g}
\newcommand{\TeichG}{\mathcal{T}_g}
\newcommand{\CirclePackGTau}{\mathsf{CPS}_{g,\tau}}
\newcommand{\CrossRatio}{{\bf c}}
\newcommand{\CrossRatioGTau}{\mathcal{C}_{g,\tau}}
\newcommand{\CrossRatioOneTau}{\mathcal{C}_{1,\tau}}
\newcommand{\DeformGTau}{\mathcal{C}_{g,\tau}}
\newcommand{\Forget}{\mathit{forg}}
\newcommand{\Uniform}{\mathit{u}}
\newcommand{\Section}{\mathit{sect}}
\newcommand{\SLTwoC}{\mathrm{SL}(2,{\mathbb C})}
\newcommand{\SLTwoR}{\mathrm{SL}(2,{\mathbb R})}
\newcommand{\SUTwo}{\mathrm{SU}(2)}
\newcommand{\PSLTwoC}{\mathrm{PSL}(2,{\mathbb C})}
\newcommand{\GLTwoZ}{\mathrm{GL}(2,{\mathbb Z})}
\newcommand{\GLTwoC}{\mathrm{GL}(2,{\mathbb C})}
\newcommand{\PSLTwoR}{\mathrm{PSL}(2,{\mathbb R})}
\newcommand{\PGLTwoR}{\mathrm{PGL}(2,{\mathbb R})}
\newcommand{\PSLTwoZ}{\mathrm{PSL}(2,{\mathbb Z})}
\newcommand{\SLTwoZ}{\mathrm{SL}(2,{\mathbb Z})}
\newcommand{\PGLTwoZ}{\mathrm{PGL}(2,{\mathbb Z})}
\newcommand{\nnn}{\noindent}
\newcommand{\C}{\mathscr C}
\newcommand{\MCG}{{\pi_0({\rm Homeo}(T))}}
\newcommand{\MMap}{{\bf \Phi}_{\mu}}
\newcommand{\HH}{{\mathbb H}^2}
\newcommand{\TT}{{\mathbb T}}
\newcommand{\X}{{\mathcal  X}}
\newcommand{\CC}{{\mathbb C}}
\newcommand{\RR}{{\mathbb R}}
\newcommand{\Q}{{\mathbb Q}}
\newcommand{\ZZ}{{\mathbb Z}}
\newcommand{\PL}{{\mathscr {PL}}}
\newcommand{\GP}{{\mathcal {GP}}}
\newcommand{\GT}{{\mathcal {GT}}}
\newcommand{\GQ}{{\mathcal {GQ}}}
\newcommand{\EE}{{{\mathcal E}(\rho)}}
\newcommand{\HHH}{{\mathbb H}^3}
\newcommand{\XBQ}{{\mathcal X}_{BQ}}

\title[The complement of the Bowditch space]{The complement of the Bowditch
space in the ${ \SLTwoC}$ character variety }

\author{Shawn Pheng Keong Ng}
\address{Department of Mathematics \\ National University of Singapore \\
2 Science Drive 2 \\ Singapore 117543}\email{shawnngpk@gmail.com}

\author{Ser Peow Tan}
\address{Department of Mathematics \\ National University of Singapore \\
2 Science Drive 2 \\ Singapore 117543} \email{mattansp@nus.edu.sg}

 \commby{}
\date{}

%\def\square{\hfill${\vcenter{\vbox{\hrule height.4pt \hbox{\vrule width.4pt
%height7pt \kern7pt \vrule width.4pt} \hrule height.4pt}}}$}
%
%\newenvironment{pf}{\noindent {\it Proof.}\quad}{\square \vskip 10pt}

% #############################################
%
%                  Abstract
%
% #############################################

\begin{abstract}
Let $\mathcal X$ be the space of type-preserving $\SLTwoC$
characters of the punctured torus $T$. The Bowditch space
${\mathcal X}_{BQ}$ is the largest open subset of ${\mathcal X}$
on which the mapping class group acts properly discontinuously,
this is characterized by two simple conditions called the
$BQ$-conditions. In this note, we show that $[\rho] \in{\rm
int}(\X \setminus \XBQ)$  if there exists an essential simple
closed curve $X$ on $T$ such that $|{\rm tr}\, \rho(X)|<0.5$.
\end{abstract}

\maketitle

\vskip 20pt
\section{{\bf Introduction}}\label{s:intro}
\vskip 20pt Let $T$ be the punctured torus and
$\pi:=\pi_1(T)=\langle X,Y\rangle$ be its fundamental group which
is free on the generators $X,Y$. The relative $\SLTwoC$ character
variety of {\em type-preserving} characters is the set
$$\X :=\{[\rho] \in{\rm Hom}(\pi, \SLTwoC)/\SLTwoC ~:~ {\rm
tr}(XYX^{-1}Y^{-1})=-2\},$$ where the equivalence is by the
conjugation action. The Bowditch space is the subset  $\XBQ
\subset \X$ of characters which satisfy two simple conditions (see
definition \ref{def:Bowditch}), this is the largest open subset of
$\X$ on which the mapping class group of $T$ acts properly
discontinuously. It is conjectured by Bowditch to be precisely the
quasi-fuchsian space ${\mathcal X}_{QF}$ (Conjecture A,
\cite{bowditch1998plms}). To attempt to verify or disprove the
conjecture, and also to study the dynamics of the action of the
mapping class group on the non-discrete characters, it is useful
to have an effective sufficient condition for $[\rho]$ to be
inside ${\rm int}(\X\setminus {\XBQ})$. We have the following:

\begin{thm}\label{thm:maintheorem}(Main theorem)
For $[\rho] \in \X$, $[\rho] \in {\rm int}(\X \setminus {\XBQ})$
if there exists $X \in \C$ such that $|{\rm tr}\,\rho(X)|<0.5$,
where $\C$ is the set of free homotopy classes of essential simple
closed curves on $T$.

\end{thm}

\begin{rmk}\label{rmk:maintheorem}~

\begin{itemize}
    \item [(a)] The bound $0.5$ in the theorem is not optimal, and can be improved, but for
    computational purposes, it is quite effective.
    \item [(b)] Jorgensen's inequality implies that if there exists $X
    \in \C$ such that  $0<|{\rm tr}\,\rho(X)|<1$, then $[\rho]$
    corresponds to a non-discrete representation.
    Rough computer experiments have shown that in fact, in
    many examples considered (no counterexamples were detected), if $|{\rm tr}\,\rho(X)|<1$
    for some $X \in \C$, with $|{\rm tr}\,\rho(X)| \not\in (-1,1)$,
    then by a trace reduction algorithm, one
    can find some $Y \in \C$ such that $|{\rm tr}\,\rho(Y)|<0.5$,
    that is, $[\rho] \in {\rm int}(\X \setminus {\XBQ})$. This
    can be regarded as supporting evidence towards Bowditch's
    conjecture as experiments with the Wada's OPTi program \cite{wadaOpti}
    has shown that in almost all cases where
    $[\rho]$ is non-discrete, there exists $X \in \C$ with $|{\rm tr}\,\rho(X)|<1$.
    \item [(c)] The theorem quantifies the result of Bowditch in
    \cite{bowditch1998plms} (Theorem 5.5) by giving an explicit bound for
    the constant $\varepsilon_0$ in his theorem, and hence generalizes Corollary 5.6 there,
     that $[\rho_0]\in {\rm int}(\X \setminus
    {\XBQ})$, where $[\rho_0]$ is the quaternionic character with
     ${\rm tr}\, \rho_0(X)={\rm tr}\, \rho_0(Y)={\rm tr}\,
     \rho_0(XY)=0$ (and hence ${\rm tr}\, \rho_0(X)=0$ for all $X \in \C$).
    \item[(d)] The set $\XBQ$ can be expected to have a very interesting and
    complicated geometry, especially at the boundary, as evidenced
    by pictures and studies of various slices of deformation spaces of
    discrete, faithful representations including the Maskit slice,
    Earle slice, Riley slice, Bers slices
    (obtained using Wada's Opti program \cite{wadaOpti}), and also the bumping phenomena
    on the boundary of the quasi-fuchsian space, as studied by various
    authors. In particular, we have the recent
    results of Bromberg that states that the closure of ${\mathcal X}_{QF}$ is
    not locally connected. Theorem \ref{thm:maintheorem}
    can be used in a computer program to draw the Bowditch space and its complement and
    this should prove useful in studying the geometry of these spaces
    and various related conjectures.
    \item[(e)] More generally, as
    studied in
    \cite{tan-wong-zhang2004necsuf}, \cite{tan-wong-zhang2004gMm} and \cite{tan-wong-zhang2004endinvariants},
    we can study the relative character
    varieties ${\X}_{\kappa}$,  where $${\rm tr}\,
    \rho(XYX^{-1}Y^{-1})=\kappa$$
    with $\kappa \neq  2$, and the Bowditch space can be defined similarly for these
    relative character varieties. If $\kappa$ is close to $-2$, our methods can be modified
    to give similar conditions for when $[\rho] \in {\rm
    int}({\X}_{\kappa}
    \setminus \XBQ)$ and this can be used together with the
    BQ-conditions to draw the Bowditch space and complement. Note
    that in this case, the Jorgensen inequality may no longer
    apply, for example if $\kappa \in (-2,2)$, since in this case
    the image may never be discrete.

\end{itemize}

\end{rmk}

%\nnn {\it Note:} This paper started out as an undergraduate
%honours project conducted by Shawn Ng under the supervision of Ser
%Peow Tan at the department of Mathematics at the National
%University of Singapore(see \cite{Ng}).

\vskip 10pt

  The rest of this paper is organized as follows. In section
  \ref{s:prelim}, we set up the notation and definitions to be used and in
  section \ref{s:proof}, we give the proof of the theorem.
   %Finally, we make some
%  concluding remarks and speculations in section \ref{s:remarks}.

\vskip 15pt

\section{{\bf Preliminaries: Notation and definitions}}\label{s:prelim}
As in the introduction, let $T$ be the punctured torus, $X,Y$ a
pair of simple closed curves on $T$ with geometric intersection
number one so that $\pi:=\pi_1(T)=\langle X,Y\rangle$. The
relative character variety of {\em type-preserving} characters is
the set (denoted by $\X$) of equivalence classes of
representations from $\pi$ to $\SLTwoC$ satisfying
\begin{equation}\label{eqn:commutator}
    {\rm tr}\,\rho (XYX^{-1}Y^{-1})=-2,
\end{equation} where two representations
are equivalent if they are conjugate by an element of $\SLTwoC$.
By classical results of Nielsen \cite{Nielsen}, (see for example
\cite{goldmanGT2003} for background and references) it does not
matter which pair of generators is used for $\pi$ in the
definition. Fixing a pair of generators $X,Y$ of $T$, by results
of Fricke, see \cite{goldmanAr2004} for an exposition, the map
\begin{equation}\label{eqn:variety}
\iota:\X \mapsto \{(x,y,z)\in {\mathbb C}^3 ~:~ x^2+y^2+z^2=xyz\},
\end{equation}
given by $$\iota[\rho]=({\rm tr}\,\rho(X),{\rm tr}\,\rho(Y),{\rm
tr}\,\rho(XY))$$ is a bijection. Henceforth we shall identify $\X$
with the cubic variety given in (\ref{eqn:variety}), and the
topology on $\X$ will be that induced by this identification. The
character $[\rho]$ such that $\iota[\rho]={\bf 0}=(0,0,0)$ is the
quaternionic character, denoted by $[\rho_0]$.

\vskip 10pt

 The outer automorphism group of $\pi$, $${\rm
Out}(\pi)={\rm Aut}(\pi)/{\rm Inn}(\pi),$$ is isomorphic to the
mapping class group of $T$ $$\pi_0({\rm Homeo}(T))\cong {\rm
GL}(2, \mathbb Z)$$  by results of Nielsen \cite{Nielsen}, and it
acts on $\X$, via the action
\begin{equation}\label{eqn:actionofOut}
    \phi([\rho])=[\rho \circ \phi^{-1}], \quad {\hbox {where}}\quad
    \phi \in {\rm Out}(\pi),~~ [\rho] \in \X.
\end{equation}
 This action is
not effective, the kernel is generated by the automorphism
$\phi_{inv}$, where $\phi_{inv}(X)=X^{-1}$,
$\phi_{inv}(Y)=Y^{-1}$, corresponding to the elliptic involution
on $T$. Denote by $\Gamma \cong {\rm PGL}(2, \mathbb Z)$ the
quotient of $\MCG$ (equivalently, ${\rm Out}(\pi)$) by the
elliptic involution, $\Gamma$ now acts effectively on $\X$.

\vskip 10pt

The set $\C$ of free homotopy classes of essential (non-trivial
and non-peripheral) simple closed curves on $T$ forms the vertices
of the pants graph $\C(T)$ of $T$, where two vertices are
connected by an edge if and only if the corresponding curves have
geometric intersection number one. $\C(T)$ is isomorphic to the
Farey graph of the hyperbolic plane, and every vertex has infinite
valence (see for example \cite{tan-wong-zhang2004gMm}). $X,Y\in
\C$ are called {\em neighbors} if they are joined by an edge in
$\C(T)$. This is equivalent to saying that $X$ and $Y$ generate
$\pi$. Note that for any $X \in \C$ and $[\rho] \in \X$, ${\rm
tr}[\rho](X)$ is well-defined. To simplify notation, we shall use
the notationally simpler ${\rm tr}\,\rho(X)$ henceforth.

\vskip 10pt

$\Gamma$ acts on $\C(T)$, and is transitive on the set of vertices
$\C$,  in fact, it is transitive on the set of neighbors $(X,Y)$,
and the set of triples of mutual neighbors $(X,Y,Z)$.

\begin{defn}\label{def:Bowditch}
The Bowditch space is the subset $\XBQ \subset \X$ consisting of
all characters $[\rho] \in \X$ satisfying the following two
conditions, called the {\em BQ-conditions}:\begin{itemize}
    \item [(i)] ${\rm tr}\,\rho(X) \not \in [-2,2]$ for any $X \in
    \C$; and
    \item [(ii)] $|{\rm tr}\,\rho(X)| \le 2$ for only finitely
    many (possibly none) $X \in \C$.
\end{itemize}
\end{defn}
In \cite{bowditch1998plms}, Bowditch showed that $\XBQ$ is open in
$\X$, and that $\Gamma$ acts properly discontinuously on $\XBQ$.
It is also not difficult to see that in fact, $\XBQ$ is the
largest open subset of $\X$ for which the action is properly
discontinuous, (see for example \cite{tan-wong-zhang2004gMm} and
\cite{tan-wong-zhang2004necsuf} for details, and generalizations
to not necessarily type-preserving characters). Furthermore, the
subset ${\X}_{QF}$ of characters corresponding to the
quasi-fuchsian representations of $\pi$ is contained in $\XBQ$ as
a connected component. Bowditch has conjectured that in fact,
${\X}_{QF}=\XBQ$.

The dynamics of the action of $\Gamma$ on ${\rm int}(\X \setminus
{\XBQ})$ is also very interesting, and some natural questions
arise. The first (see \cite{goldmanAr2005}), is whether there
exists $[\rho] \in {\rm int}(\X \setminus {\XBQ})$ such that the
closure of its orbit contains $[\rho_0]$ and intersects $\partial
\XBQ$. More generally one can ask if there is a dense orbit under
this action, or if most orbits are dense, and finally, if this
action is ergodic. Another natural question is whether ${\rm
int}(\X \setminus {\XBQ})$ is dense in $\X \setminus {\XBQ}$.

Our main theorem can be considered as a first step towards the
study of these questions as it gives an effective way of
determining if $[\rho]\in {\rm int}(\X \setminus {\XBQ})$. In
fact, the proof, which is based on a trace reduction algorithm
gives in many cases a way of constructing a sequence of elements
in the orbit of $[\rho]$ which converges to $[\rho_0]$. (In
particular, it can be modified to give an effective constant
$\varepsilon>0$ such that if there exists neighbors $(X,Y)$ such
that $|{\rm tr}\,\rho(X)|<\varepsilon$ and $|{\rm
tr}\,\rho(Y)|<\varepsilon$, then there exists a sequence of
elements in the orbit of $[\rho]$ which converges to $[\rho_0]$.
Our result is also useful for attacking the conjecture in
\cite{tan-wong-zhang2004endinvariants} that the set of ends of a
character $[\rho]$ should be a Cantor set if it contains at least
three points and is not the entire projective lamination space,
since the trace reduction algorithm given produces lots of ends of
the character when there exists $X \in \C$ with $|{\rm tr}\,
\rho(X)|<0.5$.

\section{{\bf Proof of Main Theorem: A trace reduction algorithm}}\label{s:proof}

Our proof of Theorem \ref{thm:maintheorem} is similar in spirit to
that given by Bowditch in \cite{bowditch1998plms} that $[\rho_0]
\in {\rm int}(\X \setminus {\XBQ})$, although somewhat more
geometric. The key lemma is the following:

\begin{lem}\label{lem:mainlemma}
Let $[\rho]\in \X$ and suppose that there exists $X \in \C$ such
that $|{\rm tr}\, \rho(X)|<0.5$, with ${\rm tr}\, \rho(X) \not\in
{\mathbb R}$. Then there exists a neighbor $Y$ of $X$ in $\C$ such
that $|{\rm tr}\, \rho(Y)|<|{\rm tr}\, \rho(X)|$.

\end{lem}

The theorem now follows from the lemma since if $|{\rm tr}\,
\rho(X)|<0.5$ and ${\rm tr}\, \rho(X) \in{\mathbb R}$, then
$[\rho]\not\in \XBQ$, otherwise, we can construct a sequence (of
neighbors) $\{X_n\}$ in $\C$ such that $X_0=X$,  and furthermore,
either (i) the sequence is infinite and $|{\rm tr}\,
\rho(X_{j+1})|<|{\rm tr}\, \rho(X_j)|$ for all $j$, or (ii) the
sequence is finite and terminates at $X_N$ with ${\rm tr}\,
\rho(X_{N})\in (-2,2)$. In either case, $[\rho]\not\in \XBQ$. Note
that the condition is an open condition, so $[\rho]\in {\rm
int}(\X \setminus {\XBQ})$.

\vskip 10pt

\nnn {\em Proof of Lemma \ref{lem:mainlemma}}. Let $Y_n$, $n \in
{\mathbb Z}$ denote the (successive) neighbors of $X$, and for
simplicity of notation, we use the lower case letters $x$, $y_n$
to denote ${\rm tr}\, \rho(X)$, ${\rm tr}\, \rho(Y_{n})$
respectively. The condition in the lemma is then
\begin{equation}\label{eqn:conditionforlemma}
    |x|<0.5, \qquad x
\not\in {\mathbb R}.
\end{equation}  By conjugating the representation so that
$\rho(X)$ is diagonal and $\infty$ is its attracting fixed point,
that is,
$$\rho(X)=\left(%
\begin{array}{cc}
  \lambda & 0 \\
  0 & \lambda^{-1} \\
\end{array}%
\right), \quad \rho(Y_0)=\left(%
\begin{array}{cc}
  A & B \\
  C & D \\
\end{array}%
\right),$$ we see that
\begin{equation}\label{eqn:lambda}
    x=\lambda+\lambda^{-1}, ~~{\hbox{ where}} \quad |\lambda|>1,
\end{equation}
and
\begin{equation}\label{eqn:y_n}
    y_n=A\lambda^n+D\lambda^{-n},
\end{equation}
where

\begin{equation}\label{eqn:AD=}
    AD=\frac{x^2}{x^2-4},
\end{equation} by the commutator relation (\ref{eqn:commutator}).

Write $\lambda=re^{i\theta}$, so $|\lambda|=r>1$, $\arg
\lambda=\theta \in(-\pi, \pi]$. By re-indexing and interchanging
$A$ and $D$ if necessary, we may assume that
\begin{equation}\label{eqn:D/A}
    1 \le \left |\frac{D}{A}\right | \le |\lambda|=r
\end{equation}

The idea now is that if $|x|$ is small, then $r \sim 1$ and
$|\theta| \sim \pi/2$. Hence $|A| \sim |D|\sim|x|/2$, so that
either $|y_0|<|x|$ (if $\arg A \not \sim \arg D$), or $|y_1|<|x|$
(if $\arg A \sim \arg D$). We make these arguments precise in the
following estimates.

\vskip 10pt

 From (\ref{eqn:conditionforlemma}), we have the
following bounds for $r$ and $\theta$:
\begin{equation}\label{eqn:r-bounds}
    1<r=|\lambda|<\frac{0.5+\sqrt{4.25}}{2} \cong 1.281
\end{equation}

\begin{equation}\label{eqn:theta-bound}
    -0.25<\cos \theta<0.25, \quad 0.419\pi <|\theta|<0.581\pi
    %, \qquad \cos 2\theta<-0.92
\end{equation}

From (\ref{eqn:conditionforlemma}), (\ref{eqn:AD=}) and
(\ref{eqn:D/A}), we have
\begin{equation}\label{eqn:A-bound}
    |AD| <\frac{|x|^2}{3.75} \Longrightarrow
    |A|^2<\frac{|x|^2}{3.75} \Longrightarrow
   |A|<\frac{|x|}{\sqrt{3.75}}
\end{equation}

Hence,
\begin{equation}\label{eqn:y_0}
    |y_0|=|A+D|=|A| \left |1+\frac{D}{A} \right |
    <\frac{|x|}{\sqrt{3.75}}\left |1+\frac{D}{A} \right |
\end{equation}

Now we claim that either
\begin{equation}\label{eqn:boundfor1+D/A} \left
|1+\frac{D}{A}\right |<\sqrt{3.75},\quad \hbox{or}
\end{equation}
\begin{equation}\label{eqn:boundfory_1}
    \left |\lambda+\frac{D}{A\lambda}\right |<\sqrt{3.75}
    \Longleftrightarrow
    \left|1+\frac{D\lambda^{-1}}{A\lambda}\right
    |<\frac{\sqrt{3.75}}{|\lambda|}
\end{equation}

\nnn {\it Proof of Claim:}

Suppose that the first statement is not true, that is,
$|1+\frac{D}{A} |\ge \sqrt{3.75}$. Let $D=r_1e^{i\theta_1}$,
$A=r_2e^{i\theta_2}$, and write $a:=|1+\frac{D}{A}|$,
$\alpha:=\pi-(\theta_1-\theta_2)$. So our assumption is equivalent
to
\begin{equation}\label{eqn:inequalityfora}
    a^2>3.75
\end{equation}

Applying the cosine rule to the triangle with sides corresponding
to the complex numbers $1$, $\frac{D}{A}$ and $1+\frac{D}{A}$, we
get
\begin{equation}\label{eqn:cosineruleforfirsttriangle}
    \cos (\pi-(\theta_1-\theta_2))=\cos \alpha=\frac{1+|\frac{D}{A}|^2-a^2}{2|\frac{D}{A}|}
\end{equation}
Now applying the bounds for $a$ and $|\frac{D}{A}|$ from
(\ref{eqn:inequalityfora}), (\ref{eqn:D/A}) and
(\ref{eqn:r-bounds}) to (\ref{eqn:cosineruleforfirsttriangle}) and
rounding off, we get
\begin{equation}\label{eqn:cosalpha}
    \cos \alpha <\frac{1+r^2-3.75}{2r}<-0.432
\end{equation}
%\begin{equation}\label{eqn:sinalpha}
%    |\sin \alpha|<0.866
%\end{equation}

In particular,

\begin{equation}\label{theta1theta2difference}
    |\theta_1-\theta_2|<0.36\pi, \qquad 0.64\pi<\alpha<1.36\pi.
\end{equation}

Now write $b:=|1+\frac{D\lambda^{-1}}{A\lambda}|$ and as before,
apply the cosine rule to the triangle with sides corresponding to
$1$, $\frac{D\lambda^{-1}}{A\lambda}$ and
$1+\frac{D\lambda^{-1}}{A\lambda}$ to get
\begin{equation}\label{eqn:cosineruleforsecondtriangle}
    b^2=1+\frac{1}{r^4}\left |\frac{D}{A}\right |^2-\frac{2}{r^2}\left|\frac{D}{A}\right|\cos(\alpha+2\theta)
\end{equation}
Using the bounds for $\theta$ and $\alpha$ in
(\ref{eqn:theta-bound}) and (\ref{theta1theta2difference}), we get
that
\begin{equation}\label{eqn:boundforalphaplustwotheta}
    |\alpha+2\theta|<0.522\pi \Longrightarrow
    \cos(\alpha+2\theta)>-0.07
\end{equation}
Applying (\ref{eqn:boundforalphaplustwotheta}) and  $\left
|\frac{D}{A}\right |\le r<r^2$ to
(\ref{eqn:cosineruleforsecondtriangle}), we have
\begin{equation}\label{eqn:boundforb}
    b^2<1+1+2(0.07)=2.14 ~~\Longrightarrow ~~
    b<1.463<\frac{\sqrt{3.75}}{1.281}(\cong 1.512)<\frac{\sqrt{3.75}}{|\lambda|}
\end{equation}
where the last inequality follows from (\ref{eqn:r-bounds}). This
proves the claim as the second statement of the claim holds in
this case.

\vskip 10pt

 To complete the proof of the lemma, we see that if the
first part of the claim holds, we have
$|y_0|=|A+D|<|A|\sqrt{3.75}<|x|$, otherwise,
$|y_1|=|A\lambda+D\lambda^{-1}|=|A||\lambda|b<|A|\sqrt{3.75}<|x|$,
where the last part of the inequalities in both cases follow from
(\ref{eqn:A-bound}). \qed

%\section{{\bf Concluding remarks}}\label{s:remarks}

\end{document}